%%%%%%%%%%%%%%%%%%%%%%%%%%%%%%%%%%%%%%%%%%%%%%%%%%%%%%%%%%%%%%%%%%%%%%%
%
% Modified Weak multiplier Hopf algebras.
%
% Serversie v1 van 2 juli 2014
%   Gebaseerd op de laatste vesie 1.2 van juli 2014
%
%%%%%%%%%%%%%%%%%%%%%%%%%%%%%%%%%%%%%%%%%%%%%%%%%%%%%%%%%%%%%%%%%%%%%%%

\input amstex  % Uncomment indien nodig

\input amssym
\input amssym.def

\magnification 1200
\loadmsbm
\parindent 0 cm

% Eigen standaard definities

\define\nl{\bigskip\item{}}
\define\snl{\smallskip\item{}}
\define\inspr #1{\parindent=20pt\bigskip\bf\item{#1}}
\define\iinspr #1{\parindent=27pt\bigskip\bf\item{#1}}
\define\ainspr #1{\parindent=24pt\bigskip\bf\item{#1}}
\define\aiinspr #1{\parindent=31pt\bigskip\bf\item{#1}}

\define\einspr{\parindent=0cm\bigskip}

\define\ot{\otimes}

\centerline{\bf Modified weak multiplier Hopf algebras}
\bigskip
\centerline{\it A.\ Van Daele \rm $^{(*)}$} 
\nl\nl\nl
{\bf Abstract} 
\nl 

Let $(A,\Delta)$ be a regular weak multiplier Hopf algebra (as defined and studied in [VD-W1]). Denote by $E$ the canonical idempotent of $(A,\Delta)$ and by $B$ the image $\varepsilon_s(A)$ of the source map $\varepsilon_s$. Recall that $B$ is a non-degenerate algebra, sitting nicely in the multiplier algebra $M(A)$ of $A$ so that also $M(B)$ can be viewed as a subalgebra of $M(A)$. 
\snl
Assume that $u,v$ are invertible elements in $M(B)$ so that $E(vu\ot 1)E=E$. This last condition is obviously fulfilled if $u$ and $v$ are each other inverses, but there are also other cases. Now modify $\Delta$ and define $\Delta'(a)=(u\ot 1)\Delta(a)(v\ot 1)$ for all $a\in A$. We show in this paper that $(A,\Delta')$ is again a regular weak multiplier Hopf algebra and we obtain formulas for the various data of $(A,\Delta')$ in terms of the data associated with the original pair $(A,\Delta)$.
\snl
In the case of a finite-dimensional weak Hopf algebra, the above deformation is a special case of the twists as studied in Section 6 of [N-V].
\snl
It is shown in [T-VD2] that any regular weak multiplier Hopf algebra gives rise in a natural way to a regular multiplier Hopf algebroid. This result applies to both the original weak multiplier Hopf algebra $(A,\Delta)$ and the modified version $(A,\Delta')$. However, the same method can be used to associate another regular multiplier Hopf algebroid to the triple $(A,\Delta,\Delta')$. This turns out to give an example of a regular multiplier Hopf algebroid that does not arise from a regular weak multiplier Hopf algebra although the base algebra is separable Frobenius (see [T-VD2]). 
\nl\snl
{\it July 2014} (Version 1.2)
\vskip 4 cm
\hrule
\bigskip
\parindent 0.7 cm

\item{$^{(*)}$} {\it Address}: Department of Mathematics, University of Leuven, Celestijnenlaan 200B,
B-3001 Heverlee, Belgium. {\it E-mail}: Alfons.VanDaele\@wis.kuleuven.be

\parindent 0 cm

\newpage

\bf 0. Introduction \rm
\nl
Let $A$ be a non-degenerate algebra, say over the field of complex numbers. A coproduct $\Delta$ on $A$ is a homomorphism from $A$ to the multiplier algebra $M(A\ot A)$ of the tensor product $A\ot A$ satisfying certain regularity conditions as well as coassociativity. In the setting of weak multiplier Hopf algebras, it is not assumed to be non-degenerate as this would imply the existence of a unique  extension to a homomorphism on $M(A)$ that is unital. Instead, this property is weakened in {\it the following sense}.
\snl
Consider the canonical maps $T_1$ and $T_2$ defined on $A\ot A$ by
$$T_1(a\ot b)=\Delta(a)(1\ot b)
\qquad\quad\text{and}\qquad\quad
T_2(a\ot b)=(a\ot 1)\Delta(b).$$
It is assumed that there is an idempotent $E\in M(A\ot A)$ so that the range of $T_1$ is $E(A\ot A)$ and the range of $T_2$ is $(A\ot A)E$. This idempotent is unique if it exists. The condition allows to extend the coproduct $\Delta$ to all of $M(A)$ and this extension will still be unique provided we add the requirement that $\Delta(1)=E$. In particular we have $\Delta(a)=E\Delta(a)$ and $\Delta(a)=\Delta(a)E$ for all $a\in A$. In fact, $E$ is the smallest idempotent in $M(A\ot A)$ with this property. Similarly we can extend the homomorphisms $\Delta\ot \iota$ and $\iota\ot\Delta$, where $\iota$ is the identity map, from $A\ot A$ to $M(A\ot A)$ and these extensions are again unique if we require that they map the identity of $A\ot A$ to $E\ot 1$ and $1\ot E$ respectively. Coassociativity is now expressed in the usual form $(\Delta\ot\iota)\Delta=(\iota\ot\Delta)\Delta$ and this equality even holds on elements of $M(A)$. In particular $(\Delta\ot \iota)E=(\iota\ot \Delta)E$. The {\it next requirement} is now that 
$$(\Delta\ot 1)E=(E\ot 1)(1\ot E)=(1\ot E)(E\ot 1).$$
\snl
Finally, for the pair $(A,\Delta)$ to be a weak multiplier Hopf algebra, there are also {\it conditions on the kernels} of the canonical maps $T_1$ and $T_2$. The extra assumptions are formulated entirely in terms of the canonical idempotent $E$, just as it is the case for the ranges of these maps. However, the conditions are somewhat more difficult and we refer to Definition 1.14 in [VD-W1] for a precise formulation.
\snl
A weak multiplier Hopf algebra $(A,\Delta)$ is regular if its antipode $S$ is a bijective map from $A$ to itself. In that case, also the canonical maps for the pair $(A^{\text{op}},\Delta)$ are considered (where $A^{\text{op}}$ is the algebra $A$, endowed with the opposite product). They are denoted by $T_3$ and $T_4$ and satisfy
$$T_3(a\ot b)=(1\ot b)\Delta(a)
\qquad\quad\text{and}\qquad\quad
T_4(a\ot b)=\Delta(b)(a\ot 1).$$
Also here the range of $T_3$ is $(A\ot A)E$ and the range of $T_4$ is $E(A\ot A)$. And again there as characterizations of the kernels in terms of $E$.
\snl
In this paper we only consider {\it regular weak multiplier Hopf algebras}.
\nl
For any regular weak multiplier Hopf algebra, we have the source and target maps $\varepsilon_s$ and $\varepsilon_t$ defined on $A$ by
$$\varepsilon_s(a)=\textstyle\sum_{(a)}S(a_{(1)})a_{(2)}
\qquad\quad\text{and}\qquad\quad
\varepsilon_t(a)=\textstyle\sum_{(a)}a_{(1)}S(a_{(2)})$$
where we use the Sweedler notation for the coproduct $\Delta$ on $A$ and where $S$ is the antipode of $(A,\Delta)$. The images $\varepsilon_s(A)$ and $\varepsilon_t(A)$ are non-degenerate subalgebras $B$ and $C$ of $M(A)$. They sit nicely in $M(A)$ so that also their multiplier algebras $M(B)$ and $M(C)$ can be considered as subalgebras of $M(A)$. We have that $E$ is a separability idempotent in $M(B\ot C)$ (as defined and studied in [VD3]). In fact, $B$ and $C$ can be characterized as the left and the right leg of $E$ respectively. In particular, they are commuting subalgebras in $M(A)$. 
\snl
Having recalled these various notions, we can now describe \it content of the paper\rm.
\nl
In {\it Section} 1 we start with a regular weak multiplier Hopf algebra $(A,\Delta)$. We assume that $u,v$ are invertible elements in $M(B)$ where as before, $B$ is the image of the source map. And we require that 
$$E(vu\ot 1)E=E.\tag"(0.1)"$$ 
Because $E$ is idempotent, this condition is trivially satisfied when $u$ and $v$ are each others inverses. However, as will be explained, there are also other possibilities. 
\snl
Then we define a modified coproduct $\Delta'$ on $A$ by the formula $\Delta'(a)=(u\ot 1)\Delta(a)(v\ot 1)$. This will be again a homomorphism precisely because of condition (0.1) above. It turns out that $\Delta'$ is a regular coproduct. There is still a canonical idempotent $E'$ giving the ranges of the modified canonical maps $T'_1$ and $T'_2$. It is given by $E'=(u\ot 1)E(v\ot 1)$ as expected. In this section, it is shown that the modified pair $(A,\Delta')$ is again a regular weak multiplier Hopf algebra. Also the various other data, such as the modified antipode $S'$ and the new source and target maps $\varepsilon'_s$ and $\varepsilon'_t$ are expressed in terms of the original data.  
\snl
In the case of a finite-dimensional weak Hopf algebra, the modification above turns out to be a special case of the twisting as studied by Nikshych and Vainerman in [N-V].
\snl
In {\it Section} 2 we consider the {\it associated multiplier Hopf algebroids}. In [T-VD2], it is described how any regular weak multiplier Hopf algebra gives rise to a regular multiplier Hopf algebroid. The method applies to both the original weak multiplier Hopf algebra $(A,\Delta)$ and the modified one $(A,\Delta')$. However, in this section we are interested in another multiplier Hopf algebroid coming from this example. The left multiplier bialgebroid is constructed from the original pair whereas the right multiplier bialgebroid is obtained from the new pair. The crucial property needed for this construction is the joint coassociativity rules 
$$(\Delta'\ot\iota)\Delta=(\iota\ot\Delta)\Delta'
\qquad\quad\text{and}\quad\qquad
(\Delta\ot\iota)\Delta'=(\iota\ot\Delta')\Delta.$$
In general, the regular multiplier Hopf algebroid obtained in this way, will not come from a (single) weak multiplier Hopf algebra. Nevertheless, the base algebra is separable Frobenius. This is used to answer a question raised in [T-VD2]. 
\snl
In the last section, {\it Section} 3, we draw some conclusions and we discuss possible further research on the subject.
\nl\nl
% \newpage

\it Conventions and standard references \rm
\nl
We only consider algebras over the field of complex numbers (although we believe that most of the results are still valid for more general fields). We do not require our algebras to have a unit, but we do need that the product is non-degenerate (as a bilinear form). Our algebras all turn out to be idempotent. In fact, they even have local units and this implies that the product is non-degenerate and that the algebra is idempotent.
\snl
For algebras with a non-degenerate product, it is possible to define the multiplier algebra $M(A)$. It can be characterized as the largest algebra with identity containing $A$ as an essential two-sided ideal. We use $1$ for the unit in the multiplier algebras. Of course $M(A)=A$ if and only if $A$ already has an identity. If $A$ is non-degenerate, then so is the tensor product $A\ot A$ and we consider also its multiplier algebra $M(A\ot A)$. We have natural imbeddings
$$A\ot A\subseteq M(A)\ot M(A)\subseteq M(A\ot A)$$ 
and in general, for a non-unital algebra, these two inclusions are strict. 
\snl
The opposite algebra $A^{\text{op}}$ is the algebra $A$, but endowed with the opposite product. 
\snl
We use $\iota$ for the identity map. We sometimes use $\zeta$ to denote the flip map on $A\ot A$ and for its extension to the multiplier algebra $M(A\ot A)$. A coproduct is a map from $A$ to $M(A\ot A)$. The composition of $\Delta$ with $\zeta$ is denoted by $\Delta^{\text{cop}}$. 
\snl
We will sometimes use the Sweedler notation for the coproduct. So we write $\sum_{(a)} a_{(1)} \ot a_{(2)}$ for $\Delta(a)$ when $\Delta$ is a coproduct on the algebra $A$ and $a\in A$. Some useful information about the use of the Sweedler notation in the case of non-unital algebras can be found e.g.\ in [VD2].
\snl
For the theory of weak Hopf algebras we refer to [B-N-S] and for the relation with Hopf algebroids to [B1].
\snl
For the theory of weak multiplier Hopf algebras, we refer to [VD-W1] and [VD-W2]. In particular, we find in [VD-W2] the necessary results about the source and target maps and source and target algebras that are very important for the described procedure. See also [VD-W0] for a motivational paper. The theory of multiplier Hopf algebroids has been developed recently in [T-VD1]. An important role is played by separability idempotents. For these objects, we refer to [VD3]. For the relation of weak multiplier Hopf algebras with multiplier Hopf algebroids, we refer to [T-VD2].
\snl
The conventions used in the original papers on weak multiplier Hopf algebras [VD-W1] and [VD-W2] are not always the same as in the newer paper on multiplier Hopf algebroids [T-VD1]. We will use here the conventions as in [T-VD1] and [T-VD2].
\nl\nl
\bf Acknowledgements \rm
\nl
I want to thank my colleagues and coworkers, Gabriella B\"ohm and Thomas Timmermann, for discussions about the topic of this paper. In particular, I am grateful to Gabriella B\"ohm for drawing my attention to one of her examples, related with the material in this note (see [B2]), as well as to the work of Nikshych and Vainerman on twists of finite-dimensional weak Hopf algebras (see Section 6 in [N-V]).
\nl\nl

%\newpage

%%% \input artikel1.tex % (1. Modifying a weak multiplier Hopf algebra)

\bf 1. Modifying a weak multiplier Hopf algebra \rm
\nl
We start with a regular weak multiplier Hopf algebra $(A,\Delta)$ as introduced and studied in [VD-W1]. We denote by $E$ the canonical idempotent. Its properties are used to extend the coproduct $\Delta$ to the multiplier algebra $M(A)$ and then it turns out that $E$ is actually $\Delta(1)$. There exists an invertible antipode $S$. It is an anti-automorphism of $A$ and it flips the coproduct. It gives rise to the source and target maps, given by
$$\varepsilon_s(a)=\sum_{(a)}S(a_{(1)})a_{(2)}
\qquad\quad\text{and}\quad\qquad
\varepsilon_t(a)=\sum_{(a)}a_{(1)}S(a_{(2)})$$
for $a\in A$, where we use the Sweedler notation $\Delta(a)=\sum_{(a)}a_{(1)}\ot a_{(2)}$. The range $\varepsilon_s(A)$ of the source map is denoted by $B$ and the range $\varepsilon_t(A)$ of the target map by $C$. These sets $B$ and $C$ are commuting non-degenerate subalgebras of $M(A)$. They sit nicely in $M(A)$ so that their multiplier algebras can be considered as subalgebras of $M(A)$. These multiplier algebras $M(B)$ and $M(C)$ can be characterized by
$$\align
M(B)&=\{x\in M(A) \mid \Delta(x)=E(1\ot x)\}\\
M(C)&=\{y\in M(A) \mid \Delta(y)=(y\ot 1)E\}.
\endalign$$
The algebras $B$ and $C$ are the left and right legs of $E$ respectively and because they are commuting subalgebras in $M(A)$, we will have that also $\Delta(x)=(1\ot x)E$ when $x\in M(B)$ and $\Delta(y)=E(y\ot 1)$ when $y\in M(C)$. For these results about the source and target maps and source and target algebras, we refer to [VD-W2].
\snl
Finally $E$ turns out to be a separability idempotent in $M(B\ot C)$ as studied in [VD3]. The antipodal maps $S_B:B\to C$ and $S_C:C\to B$ are characterized by the equations
$$E(x\ot 1)=E(1\ot S_B(x))
\qquad\quad\text{and}\quad\qquad
(1\ot y)E=(S_C(y)\ot 1)E$$
for $x\in B$ and $y\in C$. They are both anti-isomorphisms but need not be each others inverses. They are the restrictions of the antipode to these subalgebras (after extending the antipode from $A$ to the multiplier algebra $M(A)$). We will systematically use $S_B$ and $S_C$ for the antipode when restricted to the subalgebras $B$ and $C$. Recall that these anti-isomorphisms are completely determined by $E$ itself. Therefore, we often view $S_B$ and $S_C$ as the antipodal maps associated with the separability idempotent $E$. This remark is important as we will modify $E$ and look at the modified antipodal maps before we have modified the weak multiplier Hopf algebra. See Proposition 1.2 below.
\nl
\it Modification of the coproduct \rm
\nl
We take invertible elements $u,v\in M(B)$ and we assume that 
$$E(vu\ot 1)E=E. \tag"(1.1)"$$
This condition is of course fulfilled if $v$ and $u$ are each others inverses as $E$ is an idempotent. However it can also happen in other cases. Indeed, denote $v'=S_B(v)$ and $u'=S_C^{-1}(u)$. These elements satisfy (and are characterized by)
$$(u\ot 1)E=(1\ot u')E
\qquad\quad\text{and}\quad\qquad
E(v\ot 1)=E(1\ot v').$$
Therefore also $E(1\ot v'u')E=E$. From this we see that (1.1) is also satisfied when $u'$ and $v'$ are each other inverses.  This will not necessarily imply that also $u$ and $v$ are each others inverses.
\snl
An immediate consequence of condition (1.1) is obtained in the following proposition. In fact, it gives an equivalent formulation of this condition.

\inspr{1.1} Proposition \rm 
Using the Sweedler type notation $E=E_{(1)}\ot E_{(2)}$, we find
$$E_{(1)}vu S_C(E_{(2)})=1
\quad\quad\text{and}\quad\quad
S_B(E_{(1)})S_B(vu)E_{(2)}=1.$$

\snl \bf Proof\rm:
Because the second leg of $E$ belongs to $C$, from the formula characterizing $S_C$ we get
$$E(vu\ot 1)E=(E_{(1)}vu\ot E_{(2)})E=(E_{(1)}vu S_C(E_{(2)})\ot 1)E.$$
By assumption, this is $E$ and since $E$ is full (as in Definition 1.2 of [VD3]), it follows that $E_{(1)}vu S_C(E_{(2)})=1$. Similarly
$$E(vu\ot 1)E=E(1\ot S_B(vu))E=E(1\ot S_B(E_{(1)})S_B(vu)E_{(2)})$$
and because this is equal to $E$, again from the fullness of $E$, it follows that \newline $S_B(E_{(1)})S_B(vu)E_{(2)}=1$.
\hfill $\square$\einspr 

When $vu=1$, we recover the well-known formulas as found and explained in Proposition 1.8 of [VD3]. And just as in that case, the formulas here are given a meaning when multiplied with an element of $B$, left or right, in the first case and with an element of $C$, left or right, in the second case. Doing this, also the different steps in the proof are completely justified.
\nl
{\it We now fix such a pair} $u,v$ of invertible elements in $M(B)$ satisfying (1.1). 
\nl
We first use these elements to modify $E$.

\inspr{1.2} Proposition \rm
Define $E'=(u\ot 1)E(v\ot 1)$. Then $E'$ is a separability idempotent in $M(B\ot C)$ with associated antipodal maps $S'_B$ and $S'_C$ given by
$$S'_B(x)=S_B(vxv^{-1})
\qquad\quad\text{and}\quad\qquad
S'_C(y)=uS_C(y)u^{-1}$$
when $x\in B$ and $y\in C$.

\bf\snl Proof\rm:
i) From the assumption (1.1) we see that 
$$E'E'=(u\ot 1)E(vu\ot 1)E(v\ot 1)=(u\ot 1)E(v\ot 1)=E'$$
so that $E'$ is an idempotent. It is clear that it still belongs to $M(B\ot C)$ because $E\in M(B\ot C)$ and $u,v\in M(B)$. 
\snl
ii) We have 
$$E'(B\ot 1)=(u\ot 1)E(vB\ot 1)=(u\ot 1)E(B\ot 1)\subseteq (uB\ot C)=B\ot C.$$ 
Similarly we find $(B\ot 1)E'\subseteq B\ot C$. We also have that $(1\ot C)E'$ and $E'(1\ot C)$ are subsets of $B\ot C$. This is the first condition for $E'$ to be a separability idempotent.
\snl
iii) In a similar way, fullness of $E'$ will follow from fullness of $E$.
\snl
iv) Finally, a straightforward calculation gives
$$E'(x\ot 1)=E'(1\ot S_B(vxv^{-1}))
\quad\quad\text{and}\quad\quad
(1\ot y)E'=(uS_C(y)u^{-1}\ot 1)E'$$
for all $x\in B$ and $y\in C$. 
\snl
All of this proves that $E'$ is again a separability idempotent in $M(B\ot C)$ and that its antipodal maps $S'_B$ and $S'_C$ are given by the formulas as in the formulation of the proposition.
\hfill $\square$ \einspr

We now modify the original coproduct $\Delta$ on $A$.

\inspr{1.3} Proposition \rm
Define 
$$\Delta'(a)=(u\ot 1)\Delta(a)(v\ot 1)$$
for all $a\in A$. Then $\Delta'$ is a regular and full coproduct on $A$ with counit $\varepsilon'$ defined as 
$$\varepsilon'(a)=\varepsilon(u^{-1}av^{-1})$$
for $a\in A$ where $\varepsilon$ is the counit of the original weak multiplier Hopf algebra $(A,\Delta)$.

\snl \bf Proof\rm: 
i) It is clear that $\Delta'(a)\in M(A\ot A)$ and we have
$$\align 
\Delta'(a)\Delta'(b)
&=(u\ot 1)\Delta(a)(vu\ot 1)\Delta(b)(v\ot 1)\\
&=(u\ot 1)\Delta(a)E(vu\ot 1)E\Delta(b)(v\ot 1)\\
&=(u\ot 1)\Delta(a)E\Delta(b)(v\ot 1)\\
&=(u\ot 1)\Delta(a)\Delta(b)(v\ot 1)\\
&=(u\ot 1)\Delta(ab)(v\ot 1)=\Delta'(ab)
\endalign$$
for all $a,b\in A$.
\snl
ii) Also the four elements 
$$\align
&\Delta'(a)(1\ot b) \qquad\quad\text{and}\qquad\quad (c\ot 1)\Delta'(a) \\
&(1\ot b)\Delta'(a) \qquad\quad\text{and}\qquad\quad \Delta'(a)(c\ot 1) 
\endalign$$
belong to $A\ot A$ for all $a,b\in A$ as the same is true with $\Delta$ in the place of $\Delta'$.
\snl
iii) To prove coassociativity of $\Delta'$, take $a,b,c\in A$. On the one hand we have
$$\align
(c\ot 1\ot 1)&(\Delta'\ot\iota)(\Delta'(a)(1\ot b))\\
&=(cu\ot 1\ot 1)(\Delta\ot\iota)((u\ot 1)\Delta(a)(v\ot b))(v\ot 1\ot 1)\\
&=(cu\ot 1\ot 1)(\Delta(u)\ot 1)(\Delta\ot\iota)(\Delta(a)(1\ot b))(\Delta(v)\ot 1)(v\ot 1\ot 1).
\endalign$$
Because $u,v\in M(B)$, we have 
$$\Delta(u)=(1\ot u)E
\qquad\quad\text{and}\qquad\quad
\Delta(v)=E(1\ot v).
$$
Using these relations in the above formula, we find
$$\align
(c\ot 1\ot 1)&(\Delta'\ot\iota)(\Delta'(a)(1\ot b))\\
&=(cu\ot u\ot 1)(E\ot 1)(\Delta\ot\iota)(\Delta(a)(1\ot b))(E\ot 1)(v\ot v\ot 1)\\
&=(cu\ot u\ot 1)(\Delta\ot\iota)(\Delta(a)(1\ot b))(v\ot v\ot 1).
\endalign$$
On the other hand we find
$$\align
(\iota\ot\Delta')((c\ot 1)\Delta'(a))&(1\ot 1\ot b)\\
&=(1\ot u\ot 1)(\iota\ot\Delta)((cu\ot 1)\Delta(a)(v\ot 1))(1\ot v\ot b)\\
&=(1\ot u\ot 1)(\iota\ot\Delta)((cu\ot 1)\Delta(a))(v\ot v\ot b)\\
&=(cu\ot u\ot 1)(\Delta\ot\iota)(\Delta(a)(1\ot b))(v\ot v\ot 1)
\endalign$$
where we have used coassociativity of $\Delta$. We see that we get the same expressions and so also $\Delta'$ is coassociative.
\snl
iv) It is not hard to see that fullness of $\Delta$ gives fullness of $\Delta'$. 
\snl
v) It remains to show that $\varepsilon'$, defined as in the formulation of the proposition, is a counit for the new coproduct $\Delta'$. Clearly, we have 
$$(\varepsilon'\ot\iota)(\Delta'(a)(1\ot b))=(\varepsilon\ot\iota)(\Delta(a)(1\ot b))=ab$$
for all $a,b\in A$. Now apply $\iota\ot\varepsilon'\ot\iota$ on the equality
$$(c\ot 1\ot 1)(\Delta'\ot\iota)(\Delta'(a)(1\ot b))
=(\iota\ot\Delta')((c\ot 1)\Delta'(a))(1\ot 1\ot b).$$ 
We get 
$$(\iota\ot\varepsilon'\ot\iota)(c\ot 1\ot 1)(\Delta'\ot\iota)(\Delta'(a)(1\ot b))=(c\ot 1)\Delta'(a)(1\ot b)$$
for all $a,b,c\in A$. And because $\Delta'$ is full, we also find
$$(\iota\ot\varepsilon')((c\ot 1)\Delta'(a))=ca$$
for all $a,c$. Hence $\varepsilon'$ is the counit for the new coproduct $\Delta'$. 
\hfill$\square$\einspr

We can also define $\varepsilon''$ on $A$ by $a\mapsto\varepsilon({u'}^{-1}a{v'}^{-1})$ where as before $u',v'$ are the elements in $M(C)$ defined by
$$(u\ot 1)E=(1\ot u')E
\qquad\quad\text{and}\quad\qquad
E(v\ot 1)=E(1\ot v').$$
Because 
$$\align \Delta'(a)&=(u\ot 1)\Delta(a)(v\ot 1)\\
&=(u\ot 1)E\Delta(a)E(v\ot 1)\\
&=(1\ot u')E\Delta(a)E(1\ot v')\\
&=(1\ot u')\Delta(a)(1\ot v')
\endalign$$
for all $a$, we find 
$$(\iota\ot\varepsilon'')((c\ot 1)\Delta'(a))=(\iota\ot\varepsilon)((c\ot 1)\Delta(a))=ca$$
for all $a,c\in A$. A similar argument as above will give that $\varepsilon''$ is the counit. Hence, it has to be equal to $\varepsilon'$ as defined in the proposition.
\snl
Before we continue, we make the following remark.

\inspr{1.4} Remark \rm
i) We use this example in [T-VD2] where we study the relation  between weak multiplier Hopf algebras and multiplier Hopf algebroids. For this application, it is important to notice that there is the {\it mixed coassociativity} of the modified coproduct $\Delta'$ with respect to the original coproduct $\Delta$. This says, roughly speaking, that 
$$(\Delta'\ot \iota)\Delta=(\iota\ot\Delta)\Delta'
\qquad\quad\text{and}\qquad\quad
(\Delta\ot \iota)\Delta'=(\iota\ot\Delta')\Delta.\tag"(1.2)"$$
Strictly speaking, we have to multiply with elements of $A$ as in the case of usual coassociativity. And then the two results are easy to prove, just as the coassociativity of $\Delta'$ in Proposition 1.3.
\snl
ii) If $\Delta=\Delta'$ it follows from the uniqueness of the counit that $\varepsilon=\varepsilon'$. Also the converse is true. Indeed, assume that $\varepsilon=\varepsilon'$. Because on the one hand we have $(\iota\ot \varepsilon')\Delta'(a)=a$ while on the other hand we have $(\iota\ot \varepsilon)\Delta'(a)=uav$, we see that then $uav=a$ for all $a$. This property implies that $u$ and $v$ are each others inverses and that they belong to the center of $M(A)$. Then $\Delta'=\Delta$. This is somewhat remarkable as the equality $\varepsilon=\varepsilon'$ seems much weaker than $\Delta=\Delta'$. 
\hfill $\square$ \einspr

We will consider this further in the next section (see Theorem 2.3).
\nl
\it The ranges of the modified canonical maps \rm
\nl
We want to show that the new pair $(A,\Delta')$ is still a regular weak multiplier Hopf algebra. This involves different steps.
\snl
First recall the notations used for the four canonical maps. In the case of $(A,\Delta)$ we use
$$\align 
T_1(a\ot b)&=\Delta(a)(1\ot b) 
\qquad\qquad\qquad
T_2(a\ot b)=(a\ot 1)\Delta(b)\\
T_3(a\ot b)&=(1\ot b)\Delta(a)
\qquad\qquad\qquad
T_4(a\ot b)=\Delta(b)(a\ot 1)
\endalign$$
We will use $T'_1$, $T'_2$, $T'_3$ and $T'_4$ for the corresponding maps with $\Delta$ replaced by $\Delta'$.
\snl
We first show that $E'$ as defined in Proposition 1.2 is the canonical idempotent of $(A,\Delta')$.

\inspr{1.5} Proposition \rm
We have that 
$$\align 
\Delta'(A)(1\ot A)&=E'(A\ot A) 
\qquad\qquad\qquad
(A\ot 1)\Delta'(A)=(A\ot A)E'\\
(1\ot A)\Delta'(A)&=(A\ot A)E'
\qquad\qquad\qquad
\Delta'(A)(A\ot 1)=E(A\ot A).
\endalign$$

\snl\bf Proof\rm:
For all $a,b$ in $A$ we find
$$\Delta'(a)(1\ot b)=(u\ot 1)\Delta(a)(1\ot b)(v\ot 1)$$
and these elements span the space 
$$(u\ot 1)E(Av\ot A)=(u\ot 1)E(v\ot 1)(v^{-1}Av\ot A)=E'(A\ot A).$$
The other cases are completely similar.
\hfill $\square$ \einspr

This takes care of the necessary formulas for the ranges of the four new canonical maps. 
\snl
The next result we need is the behavior of the new coproduct $\Delta'$ on the legs of the new canonical idempotent $E'$. This is obtained in following proposition.

\inspr{1.6} Proposition \rm We have
$$(\Delta'\ot\iota)E'=(E'\ot 1)(1\ot E')=(1\ot E')(E'\ot 1).$$

\snl\bf Proof\rm:
Using the definitions of $E'$ and $\Delta'$, we find
$$\align
(\Delta'\ot\iota)E'
&=(u\ot 1\ot 1)(\Delta\ot\iota)((u\ot 1)E(v\ot 1))(v\ot 1\ot 1) \\
&=(u\ot u\ot 1)((\Delta\ot\iota)E)(v\ot v\ot 1) \\
&=(u\ot u\ot 1)((E\ot 1)(1\ot E)(v\ot v\ot 1) \\
&=(u\ot 1\ot 1)((E\ot 1)(v\ot 1\ot 1)(1\ot u\ot 1)(1\ot E)(1\ot v\ot 1) \\
&=(E'\ot 1)(1\ot E').
\endalign$$
We have used that the result holds for $E$ and $\Delta$ and also that $u$ commutes with $C$, the second leg of $E$. 
\snl
The other formula is obtained in the same way, or by using that the legs of $E'$ are the same as the legs of $E$, namely $B$ and $C$, and that these algebras commute.
\hfill $\square$ \einspr

\snl
\it The kernels of the canonical maps \rm
\nl
In Proposition 1.5 we have obtained the {\it ranges} of the four new canonical maps $T'_1$, $T'_2$, $T'_3$ and $T'_4$. In the next result, we find the kernels. It turns out that they are precisely the same as for the original maps.

\inspr{1.7} Proposition \rm
The kernels of the four modified canonical maps $T'_1$, $T'_2$, $T'_3$ and $T'_4$ coincide with the kernels of the original maps $T_1$, $T_2$, $T_3$ and $T_4$ respectively.

\snl\bf Proof\rm:
i) Take an element $\sum_i p_i\ot q_i$ in $A\ot A$. Then
$$\align T'_1(\textstyle\sum_i p_i\ot q_i)
&=\sum_i \Delta'(p_i)(1\ot q_i)\\
&=(u\ot 1)(\textstyle\sum_i \Delta(p_i)(1\ot q_i))(v\ot 1)\\
&=(u\ot 1)T_1(\textstyle\sum_i p_i\ot q_i)(v\ot 1).
\endalign$$
We see that the kernel of $T'_1$ is precisely the same as the kernel of $T_1$. The same argument will work for $T'_3$ and $T_3$. 
\snl
ii) To show that the kernel of $T'_2$ is the same as the kernel of $T_2$ we use that $\Delta'(a)$ can also be expressed as $(1\ot u')\Delta(a)(1\ot v')$. This argument will also work for $T'_4$ and $T_4$.

\hfill $\square$ \einspr

So, the kernels are left unchanged. However, the projection maps on these kernels, used in the theory to determine the new antipode, will not be the same. They will be obtained in the next proposition.
\snl
First recall that the projection maps on the kernels of the four original canonical maps are determined by 
the elements $F_1, F_2, F_3$ and $F_4$ (all in $M(A\ot A^{\text{op}})$) defined as
$$\align F_1&=(\iota\ot S)E  \qquad\quad\text{and}\qquad\quad F_3=(\iota\ot S^{-1})E \\
	F_2&=(S\ot \iota)E \qquad\quad\text{and}\qquad\quad F_4=(S^{-1}\ot \iota)E.
\endalign$$
See Proposition 4.7 in [VD-W1]. Moreover, these elements are determined by the formulas
$$\align E_{13}(F_1\ot 1) &= E_{13}(1\ot E) 
	\quad\quad\text{and}\quad\quad 
		(F_3\ot 1)E_{13}=(1\ot E)E_{13}\\
	(1\ot F_2)E_{13} &=(E\ot 1)E_{13}
	\quad\quad\text{and}\quad\quad 
	 	E_{13}(1\ot F_4)=E_{13}(E\ot 1)."
\endalign$$
See Proposition 4.6 in [VD-W1]. In the following proposition, we use the elements $F'_1, F'_2, F'_3$ and $F'_4$, defined as
$$\align F'_1&=(\iota\ot S'_C)E'  \qquad\quad\text{and}\qquad\quad F'_3=(\iota\ot {S'_B}^{-1})E' \\
	F'_2&=(S'_B\ot \iota)E' \qquad\quad\text{and}\qquad\quad F'_4=({S'_C}^{-1}\ot \iota)E'
\endalign$$
with $S'_B$ and $S'_C$ as defined in Proposition 1.2. Then we find the following formulas.

\inspr{1.8} Proposition \rm
We have
$$\align
F'_1&=(1\ot u)F_1(v\ot 1)
	\qquad\quad\text{and}\qquad\quad 
F'_2=F_2(S_B(vu)\ot 1)\\
F'_3&=(u\ot 1)F_3(1\ot v) 
	\qquad\quad\text{and}\qquad\quad 
F'_4=({S_C}^{-1}(vu)\ot 1)F_4.
\endalign$$

\snl\bf Proof\rm:
i) By definition we have $F'_1=(\iota\ot S'_C)E'$ and so, using the formula for $E'$ and the one for  $S'_C$ obtained in Proposition 1.2 we get
$$\align F'_1
&=(u\ot u)((\iota\ot S_C)E)(v\ot u^{-1})\\
&=(1\ot u)(\iota\ot S_C)((1\ot S_C^{-1}(u))E)(v\ot u^{-1})\\
&=(1\ot u)((\iota\ot S_C)E)(v\ot uu^{-1})\\
&=(1\ot u)((\iota\ot S_C)E)(v\ot 1).
\endalign$$
This proves the formula for $F'_1$.
\snl
ii) For $F'_2$, defined as $(S'_B\ot\iota)E'$ we obtain, using the formula for $S'_B$ as in Proposition 1.2, that
$$F'_2
=(S_B\ot \iota)((vu\ot 1)E(vv^{-1}\ot 1))
=((S_B\ot \iota)E)(S_B(vu)\ot 1)$$
proving the formula for $F'_2$.
\snl
iii) To obtain $F'_3$ we use that ${S'_B}^{-1}(y)=v^{-1}S_B^{-1}(y)v$ if $y\in C$. We get
$$\align 
F'_3
&=(1\ot v^{-1})(\iota\ot S_B^{-1})((u\ot 1)E(v\ot 1))(1\ot v)\\
&=(u\ot 1)(\iota\ot S_B^{-1})((E(v\ot S_B(v^{-1}))(1\ot v)\\
&=(u\ot 1)((\iota\ot S_B^{-1})E)(1\ot v)
\endalign$$
and this proves the formula for $F'_3$.
\snl
iv) Finally, to obtain $F'_4$ we use that ${S'_C}^{-1}(x)=S_C^{-1}(u^{-1}xu)$ when $x\in B$. Then we get
$$
F'_4
=({S_C}^{-1}\ot \iota)((u^{-1}u\ot 1)E(vu\ot 1))
=({S_C}^{-1}(vu)\ot 1))({S_C}^{-1}\ot \iota)E.
$$

\vskip -0.8 cm
\hfill $\square$ \einspr

With these formulas, one can verify that 
$$\align E'_{13}(F'_1\ot 1) &= E'_{13}(1\ot E') 
	\qquad\quad\text{and}\qquad\quad 
		(F'_3\ot 1)E'_{13}=(1\ot E')E'_{13}\\
	(1\ot F'_2)E'_{13} &=(E'\ot 1)E'_{13}
	\qquad\quad\text{and}\qquad\quad 
	 	E'_{13}(1\ot F'_4)=E'_{13}(E'\ot 1).
\endalign$$

\inspr{1.9} Remark \rm
There are other possible expressions for these elements $F'_i$. In the case of $F'_1$, they will all have the form $(1\ot x_1)F_1(x_2\ot 1)$ for some $x_1,x_2\in M(B)$. Similarly $F'_2$ will be of the from $(1\ot y_1)F_2(y_2\ot 1)$ where $y_1, y_2 \in M(C)$. For $F'_3$ we get expressions like $(x_1\ot 1)F_3(1\ot x_2)$with $x_1,x_2\in M(B)$ while for $F'_4$ we get them of the form $(y_1\ot 1)F_4(1\ot y_2)$ with $y_1, y_2 \in M(C)$. We will see further that this is precisely the kind of expressions we need because the kernels of the modified canonical maps are the same as the kernels of the original canonical maps.
\hfill $\square$ \einspr

We now come to the last step before we can prove the main result of this section.

\iinspr{1.10} Proposition \rm
We have the following expressions for the kernels of the modified canonical maps:
$$\text{Ker}(T'_1)=(A\ot 1)(1-F'_1)(1\ot A)
\quad\quad\quad
\text{Ker}(T'_2)=(A\ot 1)(1-F'_2)(1\ot A)\tag"(1.3)"$$
$$\text{Ker}(T'_3)=(1\ot A)(1-F'_3)(A\ot 1)
\quad\quad\quad
\text{Ker}(T'_4)=(1\ot A)(1-F'_4)(A\ot 1).\tag"(1.4)"$$

\snl \bf Proof\rm:
First remark that $\text{Ker}(T'_i)=\text{Ker}(T_i)$ for all $i$. This is shown in Proposition 1.7. On the other hand, we know that the equalities (1.3) and (1.4) hold for the original data. Therefore, we just need to argue, e.g.\ in the case of $T'_1$, that
$$(A\ot 1)(1-F_1)(1\ot A)=(A\ot 1)(1-F'_1)(1\ot A).$$
Recall that the maps $p\ot q\mapsto (p\ot 1)F_1(1\ot q)$ and $p\ot q\mapsto (p\ot 1)F'_1(1\ot q)$ are projection maps. Therefore we must show that they have the same kernel. In other words, given an element $\sum_i p_i\ot q_i$ in $A\ot A$, we must argue that
$$\sum_i(p_i\ot 1)F_1(1\ot q_i)=0$$
if and only if 
$$\sum_i(p_i\ot 1)F'_1(1\ot q_i)=0.$$
This however is true because $F'_1=(1\ot u)F_1(v\ot 1)$ as we showed in Proposition 1.8.
\snl
A similar argument works in the three other cases. Indeed, the various formulas for the maps $F'_i$ in terms of the maps $F_i$ respectively, given in Proposition 1.8, have the correct form. See also Remark 1.9 above.
\hfill$\square$\einspr

We now arrive at the main result of this section.
\nl
\it The main result \rm

\iinspr{1.11} Theorem \rm
Let $(A,\Delta)$ be a regular weak multiplier Hopf algebra with canonical idempotent $E$. Assume that $u$ and $v$ are invertible elements in $M(B)$ where $B$ is the image of the source map and that $E(vu\ot 1)E=E$. If we define $\Delta':A\to M(A\ot A)$ by $\Delta'(a)=(u\ot 1)\Delta(a)(v\ot 1)$, then $(A,\Delta')$ is again a regular multiplier Hopf algebra.

\snl\bf Proof\rm:
Consider the definition of a weak multiplier Hopf algebra as given in Definition 1.14 of [VD-W1]. Because we want to show that the modified pair $(A,\Delta')$ is regular, we also have to consider the requirements for $(A^{\text{op}},\Delta')$. This means that we must look at the four canonical maps $T'_1$, $T'_2$, $T'_3$ and $T'_4$.
\snl
We have shown in Proposition 1.3 that $\Delta'$ is a regular and full coproduct with a counit. This is the basic assumption in the definition of a regular weak multiplier Hopf algebra.
\snl
We have shown in Proposition 1.5 that we have the correct expressions for the ranges of the canonical maps as required in item i) of Definition 1.14 of [VD-W1]. In Proposition 1.6 we have obtained the required formula for the coproduct $\Delta'$ on the legs of $E'$ as formulated in item ii) of Definition 1.14 of a weak multiplier Hopf algebra in [VD-W1]. 
\snl
Finally, we obtained the necessary expressions giving the kernels of the four canonical maps in the Proposition 1.10 above. These are formulated, not in terms of the projection maps $G'_1$ and $G'_2$ as in Definition 1.14 of [VD-W1], but rather in terms of the multipliers $F'_1$ and $F'_2$ of $M(A\ot A^{\text{op}})$. This is possible as we are working in the regular case (see Section 4 of [VD-W1]). For $(A^{\text{op}},\Delta')$, we use the multipliers $F'_3$ and $F'_4$.
\hfill $\square$ \einspr

If $(A,\Delta)$ is a finite-dimensional weak Hopf algebra, then the above construction is a special case of the twists that are studied in Section 6 of [N-V]. However, one has to consider the weaker version of the condition (47) in Definition 6.1.1 of [N-V] as formulated in Remark 6.1.3.b of [N-V]. The stronger condition would imply $u=v=1$ in our case and there would be no modification at all. We also refer to Section 3 of this note for some more comments on this.
\nl
We now express the {\it data of the modified weak multiplier Hopf algebra} in terms of the original one.
\snl
We know that the original regular weak multiplier Hopf algebra $(A,\Delta)$ has an invertible antipode $S$ and also the modified one $(A,\Delta')$ has an invertible antipode $S'$. In the following proposition, we obtain the formula for the antipode $S'$ in terms of the original antipode $S$ and the modifying elements $u$ and $v$. We will give some  comments after the proof of the result and more comments at the end of this section.

\iinspr{1.12} Proposition \rm
The antipode $S'$ of the modified weak multiplier Hopf algebra $(A,\Delta')$ as obtained in Theorem 1.11 is given by
$$S'(a)=uS(vav^{-1})u^{-1}$$
for all $a$ in $A$.

\snl\bf Proof\rm:
Define $S'$ as in the formulation of the proposition. Then define $R'_1:A\ot A \to A\ot A$ by
$$R'_1(a\ot b)=\sum_{(a)}ua_{(1)}v\ot S'(a_{(2)})b$$
where we use the Sweedler notation for the original coproduct. We show that this is indeed the correct generalized inverse of $T'_1$.
\snl
i) First we show that $T'_1R'_1(a\ot b)=E'(a\ot b)$ for all $a,b\in A$. Indeed
$$\align
T'_1R'_1(a\ot b)
&=T'_1(\textstyle\sum_{(a)}ua_{(1)}v\ot S'(a_{(2)})b)\\
&=T'_1(\textstyle\sum_{(a)}ua_{(1)}v\ot uS(va_{(2)}v^{-1})u^{-1}b)\\
&=\textstyle\sum_{(a)}ua_{(1)}v\ot ua_{(2)}vuS(va_{(3)}v^{-1})u^{-1}b.
\endalign$$
We know that the property $E(vu\ot 1)E=E$ is equivalent with $E_{1}vuS_C(E_{2})=1$ (see Proposition 1.1). For any element $p\in A$ we have $\Delta(p)=\Delta(p)E$ and therefore
$$\align \sum_{(p)}p_{(1)}vuS(p_{(2)})&=\sum_{(p)}p_{(1)}E_{1}vuS(p_{(2)}E_{2})\\
&=\sum_{(p)}p_{(1)}E_{1}vuS_C(E_2)S(p_{(2)})\\
&=\sum_{(p)}p_{(1)}S(p_{(2)}).
\endalign$$
If we apply this with $vpv^{-1}$  and use that 
$$\Delta(vpv^{-1})=(1\ot v)\Delta(p)(1\ot v^{-1})$$
we will obtain
$$\sum_{(p)}p_{(1)}vuS(vp_{(2)}v^{-1})=\sum_{(p)}p_{(1)}S(vp_{(2)}v^{-1})$$
for all $p$. We now use this in the above formula with $a_{(2)}$ in the place of $p$. Then we get
$$\align
T'_1R'_1(a\ot b)
&=\textstyle\sum_{(a)}ua_{(1)}v\ot ua_{(2)}S(va_{(3)}v^{-1})u^{-1}b.\\
&=(u\ot u)(T_1R_1(vav^{-1}\ot u^{-1}b))(v\ot 1)\\
&=(u\ot u)E(vav^{-1}\ot u^{-1}b)(v\ot 1)\\
&=E'(a\ot b).
\endalign$$
For the last equality we used that $u\in M(B)$ and that it hence commutes with the second leg of $E$. And $R_1$ is the generalized inverse of $T_1$, determined by the original antipode.
\snl
ii) Next we show that $R'_1T'_1(a\ot b)=(a\ot 1)F'_1(b\ot 1)$ for all $a,b\in A$. Indeed
$$\align
R'_1T'_1(a\ot b)
&=R'_1(\textstyle\sum_{(a)}ua_{(1)}v\ot a_{(2)}b)\\
&=\textstyle\sum_{(a)}ua_{(1)}v\ot S'(ua_{(2)}v) a_{(3)}b\\
&=\textstyle\sum_{(a)}ua_{(1)}v\ot uS(vua_{(2)}vv^{-1})u^{-1} a_{(3)}b\\
&=\textstyle\sum_{(a)}ua_{(1)}v\ot uS(a_{(2)})S_B(vu)u^{-1} a_{(3)}b.
\endalign$$
Now we use that $E(vu\ot 1)E=E$ also is equivalent with 
$S_B(E_{(1)})S_B(vu)E_{(2)}=1$ (cf.\ Proposition 1.1). This in turn implies 
$$\sum_{(p)}S(p_{(1)})S_B(vu)p_{(2)}=\sum_{(p)}S(p_{(1)})p_{(2)}$$ for all $p$. Because $\Delta(u^{-1}p)=(1\ot u^{-1})\Delta(p)$  we also have 
$$\sum_{(p)}S(p_{(1)})S_B(vu)u^{-1}p_{(2)}=\sum_{(p)}S(p_{(1)})u^{-1}p_{(2)}.$$ 
We use this in the previous calculations with $p$ replaced by $a_{(2)}$. Then we get
$$\align
R'_1T'_1(a\ot b)
&=\textstyle\sum_{(a)}ua_{(1)}v\ot uS(a_{(2)})u^{-1} a_{(3)}b\\
&=(u\ot u)(R_1T_1(u^{-1}a\ot b))(v\ot 1)\\
&=(u\ot u)(u^{-1}a\ot 1)F_1(1\ot b)(v\ot 1)\\
&=(1\ot u)(a\ot 1)F_1(1\ot b)(v\ot 1)=(a\ot 1)F'_1(1\ot b).
\endalign$$
\snl
Because we have both $T'_1R'_1(a\ot b)=E'(a\ot b)$ and $R'_1T'_1(a\ot b)=(a\ot 1)F'_1(b\ot 1)$, it follows that $R'_1$ is the correct generalized inverse of $T'_1$ and hence that $S'$ is the antipode of the modified weak multiplier Hopf algebra.
\hfill $\square$ \einspr

Remark that it is a consequence of the general theory that we also must have
$$T'_2R'_2(a\ot b)=(a\ot b)E'
\qquad\quad\text{and}\qquad\quad
R'_2T'_2(a\ot b)=(a\ot 1)F'_2(1\ot b)$$
where $R'_2$ is given by 
$$R'_2(a\ot b)=\sum_{(b)}aS'(ub_{(1)}v)\ot b_{(2)}$$ 
for all $a,b\in A$. These equations can be verified using similar techniques as for $T'_1$ and $R'_1$ in the proof above.
\snl
We also see that the formula for $S'$ is compatible with the earlier formulas found for $S'_B$ and $S'_C$ in Proposition 1.2. We just have to keep in mind that $u$ and $v$ belong to $M(B)$ and that elements of $M(B)$ commute with elements of $M(C)$.
\nl
Next we obtain the formulas for the modified source and target maps.

\iinspr{1.13} Proposition \rm
The modified source and target maps $\varepsilon'_s$ and $\varepsilon'_t$ are given by
$$\varepsilon'_s(a)=u\varepsilon_s(u^{-1}a)
\qquad\quad\text{and}\qquad\quad
\varepsilon'_t(a)=\varepsilon_t(a{v'}^{-1})v'$$
for all $a\in A$ where $v'$ is the element $S_B(v)$ in $M(C)$ introduced already before.

\snl\bf Proof\rm:
i) Take $a\in A$. Then
$$\align
\varepsilon'_s(ua)
&=\textstyle\sum_{(a)}S'(ua_{(1)}v)ua_{(2)}=\textstyle\sum_{(a)}uS(vua_{(1)})u^{-1}ua_{(2)}\\
&=\textstyle\sum_{(a)}uS(vua_{(1)})a_{(2)}=\textstyle\sum_{(a)}uS(a_{(1)})S_B(uv)a_{(2)}\\
&=\textstyle\sum_{(a)}uS(a_{(1)})a_{(2)}=u\varepsilon_s(a).
\endalign$$
This proves the formula for $\varepsilon'_s$.
\snl
ii) Take again $a\in A$. Then
$$\align
\varepsilon'_t(av')
&=\textstyle\sum_{(a)}ua_{(1)}v'vS'(a_{(2)})=\textstyle\sum_{(a)}ua_{(1)}v'vuS(va_{(2)}v^{-1})u^{-1}\\
&=\textstyle\sum_{(a)}ua_{(1)}vuv'S(va_{(2)}v^{-1})u^{-1}=\textstyle\sum_{(a)}ua_{(1)}vuS(va_{(2)}v^{-1}v)u^{-1}\\
&=\textstyle\sum_{(a)}ua_{(1)}vuS(va_{(2)})u^{-1}=\textstyle\sum_{(a)}ua_{(1)}S(va_{(2)})u^{-1}=a_{(1)}S(va_{(2)})\\
&=\textstyle\sum_{(a)}a_{(1)}S(a_{(2)})S(v)=\varepsilon'_t(a)v'.
\endalign$$
This proves the formula for $\varepsilon'_t$.
\hfill $\square$ \einspr

We finish this section with a remark.
\snl
In [VD-W1] a weak multiplier Hopf algebra is defined with conditions on the ranges and the kernels of the canonical maps $T_1$ and $T_2$. From the axioms, the antipode is constructed. This follows the treatment of multiplier Hopf algebras as in [VD1]. On the other hand, just as in the case of multiplier Hopf algebras, also for weak multiplier Hopf algebras, it is possible to give a characterization in terms of a given antipode, see Theorem 2.9 in [VD-W1]. In practice, this is often more useful to verify that a given pair of an algebra with a coproduct is a weak multiplier Hopf algebra. Indeed, in many examples, there is a natural candidate for the antipode.
\snl
In the case of the modified weak multiplier Hopf algebra $(A,\Delta')$ in Theorem 1.11 however, it is not so clear from the very beginning what the antipode $S'$ should be. We can see from the proof of Proposition 1.12 that the argument is rather involved. Then the reader might be interested in the way the formula for the new antipode $S'$ was discovered. We feel it is instructive to add this here as it will allow the reader to gain a better understanding, not only of this particular case, but also of the techniques used in this theory. 
\snl
It turns out that it is easier first to find the expressions for the new source and target maps in terms of the original ones. That means that we first look at the result of Proposition 1.13 and from that result we will discover the formula in Proposition 1.12. 
\snl
Recall that the source maps $\varepsilon_s$ and $\varepsilon'_s$ satisfy
$$(\varepsilon_s\ot\iota)\Delta(a)=(1\ot a)E
\qquad\quad\text{and}\qquad\quad
(\varepsilon'_s\ot\iota)\Delta'(a)=(1\ot a)E'$$
for all $a$. This is just another way of writing $T_2R_2(b\ot a)=(b\ot a)E$ for $a,b\in A$ and then canceling $b$. And of course, similarly for $(A,\Delta')$.
\snl
From this characterization we find
$$\align (\varepsilon'_s\ot\iota)((u\ot 1)\Delta(a)(1\ot v'))
&=(1\ot a)(u\ot 1)E(1\ot v')\\
&=(u\ot 1)(\varepsilon_s\ot\iota)(\Delta(a))(1\ot v')
\endalign$$
for all $a$. If we cancel $v'$ and apply the counit, we find $\varepsilon'_s(ua)=u\varepsilon_s(a)$. This is the first formula in Proposition 1.13. 
\snl
In a completely similar way, we find the second formula.
\snl
In the next step, we try to find the formula for the antipode from the formula from the formula for $\varepsilon'_s$. We rewrite $\varepsilon'_s(ua)=u\varepsilon_s(a)$. This gives
$$\align 
\textstyle\sum_{(a)}uS(a_{(1)})a_{(2)}
&=\textstyle\sum_{(ua)}S'(u(ua)_{(1)}v)(ua)_{(2)}\\
&=\textstyle\sum_{(a)}S'(ua_{(1)}v)ua_{(2)}
\endalign$$
for all $a$. Now multiply from the right with an element $b\in A$ and replace $\textstyle\sum_{(a)}a_{(1)}\ot a_{(2)}b$ in this formula by $E_{(1)}p\ot E_{(2)}q$ where $p,q\in A$ and where as before we use the Sweedler type notation for $E$. This is justified as $\Delta(A)(1\ot A)=(A\ot A)E$. For the left hand side we get
$$uS(E_{(1)}p)E_{(2)}q=uS(p)S(E_{(1)})E_{(2)}q=uS(p)q.$$
For the right hand side we get
$$\align 
S'(uE_{(1)}pv)uE_{(2)}q
&=S'(v)S'(p)S'(uE_{(1)})E_{(2)}uq\\
&=S_B(v)S'(p)S_B(vuE_{(1)}v^{-1})E_{(2)}uq\\
&=S_B(v)S'(p)S_B(v^{-1})(S_B(E_{(1)})S_B(vu)E_{(2)})uq\\
&=S_B(v)S'(p)S_B(v^{-1})uq.
\endalign$$
In the last step, we have used  the second formula in Proposition 1.1.
We see that 
$$uS(p)q=S_B(v)S'(p)S_B(v^{-1})uq$$
and if we cancel $q$ we find the formula for the new antipode $S'$ as in Proposition 1.12. 
\snl
Using a similar argument, we find the same formula if we start from the equality $\varepsilon'_t(av')=\varepsilon_t(a)v'$, now using that $(A\ot 1)\Delta(A)=(A\ot A)E$
\snl
Observe that the formula for $\varepsilon'_s$ is obtained from the equality $T_2R_2(b\ot a)=(b\ot a)E$ while further we just used that the range of $T_1$ is $E(A\ot A)$. 
\nl
We will repeat this argument in the next section once more to find the relative antipode for the mixed case we study there (see Proposition 2.5).
\nl\nl

\bf 2. An associated multiplier Hopf algebroid \rm
\nl
Given a regular weak multiplier Hopf algebra $(A,\Delta)$, there is a natural way to associate a regular multiplier Hopf algebroid. We refer to [T-VD2] where this is proven. In the previous section, we have modified the original weak multiplier Hopf algebra by means of a pair of invertible elements $u,v\in M(B)$ satisfying the extra condition $E(vu\ot 1)E=E$. Recall that $E$ is the canonical idempotent of $(A,\Delta)$ and that $B$ is the image of the source map of $(A,\Delta)$. The new coproduct is given by $\Delta'a)=(u\ot 1)\Delta(a)(v\ot 1)$. Using the same method, we can now also associate a regular multiplier Hopf algebroid to the modified pair $(A,\Delta')$.
\snl
In this section however, we will associate a regular multiplier Hopf algebroid by using the original pair $(A,\Delta)$ for the left multiplier multiplier bialgebroid while we use the modified pair $(A,\Delta')$ for the right multiplier bialgebroid. The left and the right multiplier bialgebroids, obtained in this way, will still form a regular multiplier Hopf algebroid because we have the mixed coassociativity rules as mentioned in Remark 1.4.
\nl
We first associate a {\it left multiplier bialgebroid} to the original weak multiplier Hopf algebra. The procedure is as in Section 3 of [T-VD2]. We have the original algebra $A$ and the commuting subalgebras $B$ and $C$ in $M(A)$, together with the anti-isomorphism $S_B:B\to C$. We also consider the balanced tensor product $A\ot_\ell A$ defined by the property that 
$$xa\ot b=a\ot S_B(x)b$$ 
for all $a,b$ in $A$ and $x\in B$.

\inspr{2.1} Proposition \rm The quadruple
$$\Cal A_B:=(B,A,\iota_B, S_B)$$
is a left quantum graph. There is a left coproduct $\Delta_B$ on $A$ defined by
$$\Delta_B(a)(1\ot b)=\pi_\ell(\Delta(a)(1\ot b))$$
for all $a,b\in A$ where $\pi_\ell$ denotes the canonical projection map of $A\ot A$ onto the balanced tensor product $A\ot_\ell A$. This coproduct makes a left multiplier bialgebroid of the left quantum graph $\Cal A_B$.
\hfill $\square$ \einspr

For the proof of this result, we again refer to Section 3 of [T-VD2].
\snl
Next we associate a {\it right multiplier bialgebroid} to the modified weak multiplier Hopf algebra. We use the original algebra $A$ and the commuting subalgebras $B$ and $C$ of $M(A)$. However, now we consider the anti-isomorphism $S'_C:C\to B$ given as in Proposition 1.2 by $S'_C(y)=uS_C(y)u^{-1}$ for $y\in C$. The relevant balanced tensor product, denoted by $A\ot'_r A$, is defined by the property that 
$$a\ot by=aS'_C(y)\ot b$$ 
whenever $a,b\in A$ and $y\in C$.

\inspr{2.2} Proposition \rm The quadruple 
$$\Cal A'_C:=(C,A,\iota_C,S'_C)$$ 
is a right quantum graph. There is a right coproduct $\Delta'_C$ on $A$ defined by 
$$(c\ot 1)\Delta'_C(a)=\pi'_r((c\ot 1)\Delta'(a))$$
for all $a,c\in A$ where $\pi'_r$ denotes the canonical projection map of $A\ot A$ onto the balanced tensor product $A\ot'_r A$. This coproduct makes a right multiplier bialgebroid of the right quantum graph $\Cal A'_C$.
\hfill $\square$ \einspr

Again this result follows from the general theory as developed in [T-VD2].
\snl
We will now show that the left and the right multiplier bialgebroids together form a regular multiplier Hopf algebroid as in Definition 6.4 of [T-VD1].

\inspr{2.3} Theorem \rm
The left and the right quantum graphs $\Cal A_B$ and $\Cal A'_C$, as given in Proposition 2.1 and Proposition 2.2, are compatible. The left multiplier bialgebroid $(\Cal A_B, \Delta_B)$ together with the right multiplier bialgebroid $(\Cal A'_C,\Delta'_C)$ form a regular multiplier Hopf algebroid.

\snl\bf Proof\rm:
We know that all the canonical maps under consideration are bijective. So we just have to prove that the two forms of joint coassociativity hold. 
\snl
We first consider
$$(\Delta_B\ot\iota)((1\ot b)\Delta'_C(a))(c\ot 1\ot 1)
=(1\ot 1\ot b)(\iota\ot \Delta'_C)(\Delta_B(a)(c\ot 1))\tag"(2.1)"$$
for all $a,b,c\in A$. 
\snl
For the left hand side of (2.1) we have
$$\align 
(\Delta_B\ot\iota)((1\ot b)\Delta'_C(a))&(c\ot 1\ot 1)\\
&=(\pi_\ell\ot\iota)((\Delta\ot\iota)(\pi'_r((1\ot b)\Delta'(a)))(c\ot 1\ot 1)).
\endalign$$
We claim that (in an appropriate sense)
$$(\Delta\ot\iota)\circ\pi'_r=(\iota\ot\pi'_r)\circ(\Delta\ot\iota).\tag"(2.2)"$$
Then the left hand side of (2.1) will be
$(\pi_\ell\ot \iota)(\iota\ot\pi'_r)$
applied to 
$$(\Delta\ot\iota)((1\ot b)\Delta'(a))(c\ot 1\ot 1).$$
Similarly, the right hand side will turn out to be the image of
$$(1\ot 1\ot b)((\iota\ot \Delta')(\Delta(a)(c\ot 1))$$
under the composition $(\iota\ot\pi'_r)(\pi_\ell\ot \iota)$.
Now we have seen in Remark 1.4 that
$$(\Delta\ot\iota)((1\ot b)\Delta'(a))(c\ot 1\ot 1)
=(1\ot 1\ot b)(\iota\ot \Delta')(\Delta(a)(c\ot 1)).$$
Moreover, it is easy to see that the maps $\pi_\ell\ot \iota$ and $\iota\ot\pi'_r$ commute. Indeed, when we look at the defining properties, we see that in the first case we have the relations $xa\ot b\ot c=a\ot S_B(x)b\ot c$ for $x\in B$ and in the second case $a\ot b\ot cy=a\ot bS'_C(y)\ot c$ for $y\in C$. So we just need to observe that left multiplication by an element of $C$ and right multiplication by an element of $B$ commute.
\snl
So, in order to complete the proof of (2.1), we just need an argument for (2.2). To show this consider $p\ot q\in A\ot A$. On the one hand we have 
$$\pi'_r(p\ot qy))=\pi'_r(pS'_C(y)\ot q)$$
when $y\in C$. On the other hand we then get 
$$(\Delta\ot\iota)(pS'_C(y)\ot q)=((\Delta\ot\iota)(p\ot q))((1\ot S'_C(y)\ot 1)$$
because $S'_C(y)\in B$. This will precisely mean that $\Delta\ot\iota$ is well defined on $\pi'_r(A\ot A)$ and that (2.2) holds.
\snl
This completes the proof of (2.1). 
\snl
In a completely similar way, one can show that also
$$(c\ot 1\ot 1)(\Delta'_C\ot\iota)(\Delta_B(a)(1\ot b))
=((\iota\ot \Delta_B)((c\ot 1)\Delta'_C(a)))(1\ot 1\ot b)$$
for all $a,b,c\in A$. Then we have proven the theorem.
\hfill $\square$ \einspr

We know from the general theory that there is an antipode for this regular multiplier Hopf algebroid. We will denote it by $S_r$ while we keep $S$ for the original antipode of the given weak multiplier Hopf algebra $(A,\Delta)$ and $S'$ for the antipode of the modified weak multiplier Hopf algebra $(A,\Delta')$ as obtained in Proposition 1.12 in the previous section. We think of $S_r$ as the {\it relative} antipode.
\snl
We recall some of the properties of the antipode $S_r$ of the regular multiplier Hopf algebroid (see Definition 6.6  and Theorem 6.8 in [T-VD1]) in the next proposition. 

\inspr{2.4} Proposition \rm
The antipode $S_r$ is an anti-isomorphism of $A$, it coincides with $S_B$ on $B$ and with $S'_C$ on $C$. Moreover it satisfies
$$T_\rho(\iota\ot S_r)_\rho T=\iota\ot S_r
\qquad\quad\text{and}\qquad\quad
_\lambda T(S_r\ot\iota)T_\lambda=S_r\ot\iota.\tag"(2.3)"$$

\vskip -0.8 cm
\hfill $\square$ \einspr

Recall the formulas for these canonical maps, with the notations as in [T-VD1]:
$$\align 
&T_\lambda:a\ot b\mapsto \Delta_B(b)(a\ot 1)
\qquad\quad\quad\qquad
T_\rho: a\ot b\mapsto \Delta_B(a)(1\ot b)\\
&_\lambda T: a\ot b \mapsto (a\ot 1)\Delta'_C(b)
\qquad\quad\quad\qquad
_\rho T:a\ot b \mapsto (1\ot b)\Delta'_C(a)
\endalign$$
for all $a,b$.
\snl
We need to compare these maps with the canonical maps, defined from $A\ot A$ to itself, associated with the two regular weak multiplier Hopf algebras $(A,\Delta)$ and $(A,\Delta')$ with the notations as in [VD-W1]:
$$\align 
&T_1:a\ot b\mapsto\Delta(a)(1\ot b)
\qquad\quad\quad\qquad
T'_2: a\ot b \mapsto (a\ot 1)\Delta'(b)\\
&T'_3: a\ot b \mapsto (1\ot b)\Delta'(a)
\qquad\quad\quad\qquad
T_4: a\ot b \mapsto \Delta(b)(a\ot 1)
\endalign$$
for all $a,b$. We have the correspondence of the pair $(T_1,T_4)$ with $(T_\rho,T_\lambda)$ on the one hand and of the pair $(T'_2,T'_3)$ with $(_\lambda T,_\rho T)$ on the other hand.

\inspr{2.5} Proposition \rm 
The antipode $S_r$ of the regular multiplier Hopf algebroid in Theorem 2.3 is given by
$S_r(a)=uS(a)u^{-1}$ for all $a$, where $S$ is the antipode of the original weak multiplier Hopf algebra $(A,\Delta)$

\snl\bf Proof\rm:
We know that the antipode $S_r$ satisfies the equality
$$_\lambda T(S_r\ot\iota)T_\lambda=S_r\ot\iota.$$
Recall that the various maps in the above equality are all bijections between appropriate balanced tensor products. 
Therefore, this property turns out to be equivalent with the equation
$$T'_2(S_r\ot\iota)T_4(a\ot b)=(S_r(a)\ot b)E'\tag"(2.4)"$$
of maps from $A\ot A$ to itself, where $a,b\in A$. This is true because the range of $T'_2$ is $(A\ot A)E'$.
For the left hand side of (2.4) we find
$$\align
T'_2(S_r\ot\iota)T_4(a\ot b)
&=T'_2(S_r\ot\iota)(\textstyle\sum_{(b)}b_{(1)}a\ot b_{(2)})\\
&=T'_2(\textstyle\sum_{(b)}S_r(b_{(1)}a)\ot b_{(2)})\\
&=\textstyle\sum_{(b)}(S_r(b_{(1)}a)\ot 1)\Delta'(b_{(2)})\\
&=\textstyle\sum_{(b)}S_r(b_{(1)}a)ub_{(2)}v\ot b_{(3)}.
\endalign$$
And for the right hand side of (2.4) we obtain
$$\align
(S_r(a)\ot b)E'
&=(S_r(a)u\ot b)E(v\ot 1) \\
&=\textstyle\sum_{(b)}S_r(a)uS(b_{(1)})b_{(2)}v\ot b_{(3)}
\endalign$$
If we apply the counit on the second factor, the equality (2.4) now yields
$$\textstyle\sum_{(b)}S_r(b_{(1)}a)ub_{(2)}v=
\textstyle\sum_{(b)}S_r(a)uS(b_{(1)})b_{(2)}v$$
for all $a,b\in A$. Now, we cancel $v$ and we use that $S_r$ is an anti-isomorphism and we write the left hand side as
$$\textstyle\sum_{(b)}S_r(a)S_r(b_{(1)})ub_{(2)}.$$  
Because the equation is true for all $a$, we get 
$$\textstyle\sum_{(b)}S_r(b_{(1)})ub_{(2)}=
\textstyle\sum_{(b)}uS(b_{(1)})b_{(2)}$$
for all $b$. We multiply with $c$ from the right and we replace
replace $\Delta(b)(1\ot c)$ by $E(p\ot q)$. We find that
$$S_r(E_{(1)}p)uE_{(2)}q=uS(E_{(1)}p)E_{(2)}q$$
for all $p,q\in A$. We can cancel $q$. We use further that $S$ an $S_r$ are  anti-isomorphisms. We get
$$S_r(p)S_r(E_{(1)})uE_{(2)}=uS(p)S(E_{(1)})E_{(2)}=uS(p)$$
for all $p$. With $p=1$ we get $S_r(E_{(1)})uE_{(2)}=u$ and if we insert this, we finally find $S_r(p)u=uS(p)$.
\snl
This completes the proof.
\hfill $\square$ \einspr

Compare this argument with the one we gave at the end of the previous section where it is explained how we found the modified antipode $S'$.
\snl
One verifies immediately that $S_r$ coincides with $S_B$ on $B$ and with $S'_C$ on $C$. Indeed, if $x\in B$ then 
$$S_r(x)=uS_B(x)u^{-1}=S_B(x)$$ 
because $S_B(x)\in C$ and $u\in M(B)$ and hence these elements commute. On the other hand, if $y\in C$ we have  
$$S_r(y)=uS(y)u^{-1}=uS(vyv^{-1})u^{-1}=S'(y)=S'_C(y)$$ 
using that  $y\in C$ and $v\in M(B)$ and so these elements also commute. We see that $S_r(y)=S'_C(y)$. 
\snl
From the general theory, we know that the antipode $S_r$ also has to satisfy the other equation formulated in Proposition 2.4, namely
$$T_\rho(\iota\ot S_r)_\rho T=\iota\ot S_r.$$ 
Also here we have bijections between various balanced tensor products. 
Again this equality is equivalent with 
$$T_1(\iota\ot S_r)T'_3(a\ot b)=E(a\ot S_r(b)),$$ 
now in $A\ot A$, for all $a,b$. The reason is that the range of $T_1$ is $E(A\ot A)$. 
To verify this equation, take $a,b\in A$. Then 
$$\align 
T_1(\iota\ot S_r)T'_3(a\ot b)
&=T_1(\iota\ot S_r)(\textstyle\sum_{(a)}ua_{(1)}v\ot ba_{(2)})\\
&=T_1(\textstyle\sum_{(a)}ua_{(1)}v\ot S_r(ba_{(2)}))\\
&=\textstyle\sum_{(a)}\Delta(ua_{(1)}v)(1\ot S_r(ba_{(2)}))\\
&=\textstyle\sum_{(a)}a_{(1)}\ot ua_{(2)}vS_r(ba_{(3)})\\
&=\textstyle\sum_{(a)}a_{(1)}\ot ua_{(2)}vuS(ba_{(3)})u^{-1}\\
&=\textstyle\sum_{(a)}a_{(1)}\ot ua_{(2)}vuS(a_{(3)})S(b)u^{-1}\\
&=\textstyle\sum_{(a)}a_{(1)}\ot ua_{(2)}S(a_{(3)})S(b)u^{-1}\\
&=\textstyle\sum_{(a)}a_{(1)}\ot a_{(2)}S(a_{(3)})uS(b)u^{-1}\\
&=E(a\ot S_r(b)).
\endalign$$

From the formula for the antipode $S_r$, given in Proposition 2.5, we can obtain the formulas for the counits $\varepsilon_B$ and $\varepsilon'_C$ of the multiplier Hopf algebroid in Theorem 2.3. In accordance with previous conventions, we use $\varepsilon_C$ for the right counit of the multiplier Hopf algebroid associated with the original weak multiplier Hopf algebra $(A,\Delta)$.
\snl
We know that the left and right counits $\varepsilon_B$ and $\varepsilon_C$ of the multiplier Hopf algebroid associated with $(A,\Delta)$ are given by 
$$\varepsilon_B(a)=S^{-1}(\varepsilon_t(a))
\qquad\quad\text{and}\qquad\quad
\varepsilon_C(a)=S^{-1}(\varepsilon_s(a))$$
for all $a\in A$ where $\varepsilon_s$ and $\varepsilon_t$ are the source and target maps of the original weak multiplier Hopf algebra $(A,\Delta)$. 
See e.g.\ Proposition 3.6 and Proposition 3.7 in [T-VD2].

\inspr{2.6} Proposition \rm
The right counit $\varepsilon'_C$ of the multiplier Hopf algebroid in Theorem 2.3 satisfies
$\varepsilon'_C(a)=\varepsilon_C(u^{-1}au)$ for all $a$.

\snl \bf Proof\rm:
For the right counital map  $\varepsilon'_C$, as defined in Definition 6.6 of [T-VD1], we have here
$$S'_C(\varepsilon'_C(a))=\textstyle\sum_{(a)}S_r(a_{(1)})a_{(2)}=\textstyle\sum_{(a)}uS(a_{(1)})u^{-1}a_{(2)}\tag"(2.1)"$$
for all $a$. As before, $S$ is the antipode of the original weak multiplier Hopf algebra and the Sweedler notation is also used for the original coproduct.
\snl
Recall from Proposition 1.2 that $S'_C(y)=uS_C(y)u^{-1}$ when $y\in C$. This implies that ${S'_C}^{-1}(x)={S_C}^{-1}(u^{-1}xu)$ for $x\in B$. If we use this formula in 
(2.1) above, we find 
$$\varepsilon'_C(a)=\textstyle\sum_{(a)}{S_C}^{-1}(u^{-1}uS(a_{(1)})u^{-1}a_{(2)}u)=S^{-1}(\varepsilon_s(u^{-1}au))$$
and this is precisely $\varepsilon_C(u^{-1}au)$.
\hfill $\square$ \einspr

It is also instructive to verify the formula for $\varepsilon_B$ as given in Definition 6.6 in [T-VD1]. In this case, it says
$$\align 
S_B(\varepsilon_B(a))&=\textstyle\sum_{(a)}ua_{(1)}vS_r(a_{(2)})\\
&=\textstyle\sum_{(a)}ua_{(1)}vuS(a_{(2)})u^{-1}\\
&=\textstyle\sum_{(a)}ua_{(1)}S(a_{(2)})u^{-1}\\
&=u\varepsilon_t(a)u^{-1}=\varepsilon_t(a)
\endalign$$
for all $a$. This indeed gives $\varepsilon_B(a)=S^{-1}(\varepsilon_t(a))$ for all $a\in A$.
\snl
We can also verify the formula $S_r(\varepsilon_B(a))=\varepsilon'_C(S_r(a))$ for $a$ in $A$. Indeed, for the left hand side we have
$$S_r(\varepsilon_B(a))=uS(\varepsilon_B(a))u^{-1}=u\varepsilon_C(S(a))u^{-1}=\varepsilon_C(S(a)),$$
while for the right hand side we have
$$\varepsilon'_C(S_r(a))=\varepsilon'_C(uS(a)u^{-1})=\varepsilon_C(u^{-1}uS(a)u^{-1}u)=\varepsilon_C(S(a)).$$

Let us finish this section with another remark.
\snl
In Section 1, at the end, we explained how the antipode $S'$ for the modified weak multiplier Hopf algebra was found. But we did not use this argument to actually prove that the candidate for $S'$ was the correct one. In this section, we used another method. We really obtained the required formula for the new antipode from its properties. We could also have started with the formula and showed that it satisfies the requirements of the antipode. This would be more involved as the defining properties of the antipode in the case of a multiplier Hopf algebroid are more complicated than in the case of a weak multiplier Hopf algebra. Using the method as we did here uses already the general result in the theory of multiplier Hopf algebroids that the antipode exists.
\nl\nl

\bf 3. Conclusions and further research \rm
\nl
In this paper, we constructed a new regular weak multiplier Hopf algebra by modifying the coproduct of a given regular weak multiplier Hopf algebra $(A,\Delta)$. The new coproduct $\Delta'$ is given in terms of the original coproduct by
$$\Delta'(a)=(u\ot 1)\Delta(a)(v\ot 1)$$
for all $a\in A$. The elements $u,v$ are invertible elements in the multiplier algebra $M(B)$ of the base algebra $B$ satisfying 
$$E(vu\ot 1)E=E.\tag"(3.1)"$$
Recall that $E$ is the canonical idempotent of $(A,\Delta)$ and that $B$ is the image of the source map. The algebra $B$ can also be characterized as the left leg of $E$ (with the right interpretation). The data of the modified pair $(A,\Delta')$ are given in terms of the data of the original pair.
\snl
As we mentioned already in Section 1, in a remark made after Theorem 1.11, if the pair $(A,\Delta)$ is a finite-dimensional weak Hopf algebra, the modified pair we consider is a special case of the twists as studied in Section 6 of [N-V], taking into account the remark 6.1.3.b of that paper. This suggest the possibility of generalizing the results of this note by allowing more general twists of the weak multiplier Hopf algebra.
\snl
It will be interesting to see how this modifying procedure (and its possible generalization) affects the possible existence of integrals and the dual construction. The modified weak multiplier Hopf algebra will have a faithful set of integrals if that is the case for the original one. Then it will be possible to consider the dual of the modified pair. It will be some kind of modification of the dual, but of another type. It is expected that this type of modification will also be possible without the existence of integrals. This should be investigated.
\snl 
Our construction provides new examples of regular weak multiplier Hopf algebras. It is worthwhile to look for concrete examples. Also for the more general construction and its dual version above.
\snl
There seems to be a {\it more general type} of modification. For this one takes any $U,V$ in $M(A\ot A)$, but with the appropriate conditions so that the new coproduct $\Delta'$ defined on $A$ by $\Delta'(a)=U\Delta(a)V$ still makes $A$ into a (regular) weak multiplier Hopf algebra. See Section 6.1 in [N-V].
\snl
There are also some interesting {\it special cases} that should be investigated.
\snl
There is the case of a weak multiplier Hopf $^*$-algebra. In that case, we want $u$ to be unitary and $v=u^*$ in order to have that $\Delta'$ is still a $^*$-homomorphism. 
\snl
Another special case is when $B$ is abelian, or more generally, when $u$ and $v$ belong to the center of $M(B)$. In that case, the condition (3.1) forces $u$ and $v$ to be each others inverses. This implies that the new canonical idempotent $E'$ will coincide with the original one. It also follows that the antipodal maps $S_B$ and $S_C$ remain the same. However, the element $u$ need not be in the center of $A$ and so the modified coproduct can still be different form the original one. An interesting feature of this case is also that all the modified canonical maps have the same ranges and kernels as the original ones, with the same choices of projections maps on these sets, but still can be different maps.
\snl
We also have considered the special multiplier Hopf algebroid associated with the triple $(A,\Delta,\Delta')$. Also here more concrete examples should be constructed and special cases can be further investigated, just as above for the modified weak multiplier Hopf algebra. 
\snl
This situation also raises {\it two questions}. 
\snl
There is first a question about the general theory of multiplier Hopf algebroids. Is it possible to impose a condition that guarantees that the regular multiplier Hopf algebroid with a base algebra that is separable Frobenius, will necessarily come from a regular weak multiplier Hopf algebra?  It should be a condition that is formulated in terms of the data of the multiplier Hopf algebroid, without reference to the property of the base algebra being separable Frobenius. 
\snl
A second question concerns the general theory of weak multiplier Hopf algebras. Is it possible to develop a more general theory where the algebra $A$ is considered with two coproducts $\Delta_1,\Delta_2$, satisfying mixed coassociativity
$$(\Delta_1\ot\iota)\Delta_2=(\iota\ot\Delta_2)\Delta_1
\qquad\quad\text{and}\qquad\quad
(\Delta_2\ot\iota)\Delta_1=(\iota\ot\Delta_1)\Delta_2$$
and with an antipode $S_r$ satisfying 
$$\Delta_1(S_r(a))=\zeta(S_r\ot S_r)\Delta_2(a)
\quad\quad\text{and}\quad\quad
\Delta_2(S_r(a))=\zeta(S_r\ot S_r)\Delta_1(a)$$
for all $a$? Recall that $\zeta$ is the flip map. Any regular multiplier Hopf algebroid, with a separable Frobenius base algebra should come from such a  generalized regular weak multiplier Hopf algebra $(A,\Delta_1,\Delta_2)$.
\nl\nl

\bf References \rm
\nl
{[\bf B1]} G.\ B\"ohm: {\it An alternative notion of Hopf algebroid}. In "Hopf algebras in noncommutative geometry and physics",
pp. 31-53, Lecture Notes in Pure and Appl. Math., 239, Dekker, New York, 2005. See also arXiv: 0301.169 [math.QA].
\snl
{[\bf B2]} G.\ B\"ohm: {\it Hopf algebroids}.  Preprint version of a chapter for Handbook of Algebra, arXiv: 0805.3806 [math.QA].
\snl
{[\bf B-N-S]} G.\ B\"ohm, F.\ Nill \& K.\ Szlach\'anyi: {\it Weak Hopf algebras I. Integral theory and C$^*$-structure}. J.\ Algebra 221 (1999), 385-438. 
\snl
{[\bf N-V]} D.\ Nikshych \& L.\ Vainerman: {\it Finite quantum groupoids and their applications}. In {\it New Directions in Hopf algebras}. MSRI Publications, Vol.\ 43 (2002), 211-262. See also arXiv: 0006.057v2 [math.QA].
\snl
{[\bf VD1]} A.\ Van Daele: {\it Multiplier Hopf algebras}. Trans. Am. Math. Soc.  342(2) (1994), 917-932.
\snl
{[\bf VD2]} A.\ Van Daele: {\it Tools for working with multiplier Hopf algebras}. ASJE (The Arabian Journal for Science and Engineering) C - Theme-Issue 33 (2008), 505--528.  See also arXiv: 0806.2089 [math.QA]. 
\snl
{[\bf VD3]} A.\ Van Daele: {\it Separability idempotents and multiplier algebras}. Preprint University of Leuven (2013), arXiv: 1301.4398v1 [math.RA].
\snl
{[\bf VD-W0]} A.\ Van Daele \& S.\ Wang: {\it Weak multiplier Hopf algebras. Preliminaries and basic examples}. Operator algebras and quantum groups. Banach Center Publication 98 (2012), 367-415. See also arXiv: 1210.3954v1 [math.RA].
\snl
{[\bf VD-W1]} A.\ Van Daele \& S.\ Wang: {\it Weak multiplier Hopf algebras I. The main theory}. Preprint University of Leuven and Southeast University of Nanjing (2012), arXiv: 1210.4395v1 [math.RA]. Journal für die reine und angewandte Mathematik (Crelles Journal), to appear.
\snl
{[\bf VD-W2]} A.\ Van Daele \& S.\ Wang: {\it Weak multiplier Hopf algebras II. The source and target algebras}. Preprint University of Leuven and Southeast University of Nanjing (2014), arXiv: 1403.7906v1 [math.RA].
\snl
{[\bf T-VD1]} T.\ Timmermann \& A.\ Van Daele: {\it Regular multiplier Hopf algebroids. Basic theory and examples}. Preprint University of M\"unster and University of Leuven (2013), arXiv: 1307.0769v3 [math.QA].
\snl
{[\bf T-VD2]} T.\ Timmermann \& A.\ Van Daele: {\it Multiplier Hopf algebroids arising from weak multiplier Hopf algebras}. Preprint University of M\"unster and University of Leuven (2014), arXiv: 1406.3509v1 [math.RA].

\end